\documentclass[reqno]{article}
\usepackage{amssymb, amsmath, amsthm, amsfonts, amscd}
\usepackage[english]{babel}

\usepackage{graphicx}

\usepackage{epstopdf}

\usepackage{breqn}
\usepackage{mathtools}
\usepackage{color}
\usepackage{hyperref}
\hypersetup{
    colorlinks=true,
    citecolor=red,
    filecolor=black,
    linkcolor=blue,
    urlcolor=black
}

\chardef\bslash=`\\ 

\newtheorem{theorem}{Theorem}[section]
\newtheorem{corollary}[theorem]{Corollary}
\newtheorem{lemma}[theorem]{Lemma}
\newtheorem{proposition}[theorem]{Proposition}

\newtheorem{thmx}{Theorem}

\newtheorem{question}{Question}[section]

\newcommand{\N}{\mathbb{N}}

\newcommand{\Z}{\mathbb{Z}}

\newcommand{\R}{\mathbb{R}}

\newcommand{\T}{\mathbb{T}}

\def\a{\alpha }
\def\g{\gamma }

\def\l{\lambda }

\def\r{\rho}
\def\s{\sigma}
\def\t{\tau}

\def\w{\omega}
\def\sm{C^{\infty} }
\def\e{\varepsilon}
\def\f{\varphi}
\def\.{\cdot }
\def\ra{\rightarrow}

\def\begeq{\begin{equation*}}
\def\endeq{\end{equation*}}

\title{
\textsc{\textbf{Local Rigidity of Diophantine translations in higher dimensional tori}}\\
\author{Nikolaos Karaliolios \footnote{Imperial College London. Email: n.karaliolios@imperial.ac.uk}
}
}

\begin{document}

\maketitle

\begin{abstract}
We prove a theorem asserting that, given a Diophantine
rotation $\alpha $ in a torus $\T ^{d} \equiv \R ^{d} / \Z ^{d}$,
any perturbation, small enough in the $C^{\infty}$ topology,
that does not
destroy all orbits with rotation vector $\alpha$ is actually
smoothly conjugate to the rigid rotation. The proof relies
on a K.A.M. scheme (named after Kolmogorov-Arnol'd-Moser),
where at each step the existence of an invariant measure with rotation
vector $\alpha$ assures that we can linearize the equations
around the same rotation $\alpha$. The proof of the convergence of
the scheme is carried out in the $C^{\infty}$ category.
%
%
\end{abstract}

\tableofcontents

\section{Introduction}

Let $R_{\a } : x \mapsto x+\a \mod \Z ^{d}$ be a translation
on a torus $\T ^{d}$ with $d \in \N^{*}$. The search for conditions
under which a diffeomorphism $f \in \mathrm{Diff}^{\infty}(\T ^{d})$
is guaranteed to be smoothly conjugate to $R_{\a}$ is a
very old subject in dynamical systems and the source
of very deep and far-reaching studies, see for example
\cite{Herm79} and \cite{YocAst} for the case $d=1$.

To our best knowledge, the strongest rigidity result  on
perturbations of Diophantine rotations in higher dimensional
tori in the literature is the one proved in \cite{Herm79}.
This theorem, apart from the smallness assumptions, needs
the preservation of a volume form, something that assures
that \textit{every} orbit rotates at the speed of the
Diophantine rotation, so that the analogy with the
one-dimensional theory is direct.

Our goal in the present article is to relax the condition of
preservation of a (harmonic) volume form to a considerably
weaker one, which seems partly optimal. The closeness-to-rotations
condition is, a priori at least, not indispensable, while the
Diophantine property is known to be thus, since in the
Liouvillean world rigid rotations tend to be fragile. The present
rigidity theorem, whose precise statement is given in thm
\ref{thm rig}, is in fact an instance of the strength of
the Diophantine condition and of the K.A.M. machinery.
\begin{thmx} \label{thm A}
Let $\a \in \T ^{d}$, $d \in \N ^{*}$, be a Diophantine rotation and
$f \in \mathrm{Diff}^{\infty} (\T ^{d}) $ be a small enough
perturbation. Then, if $\a $ is in the convex hull of the
rotation set of $f$, the diffeomorphism $f$ is
smoothly conjugate to the translation by $\a $.
\end{thmx}

The motivation for this theorem comes from a conjecture
concerning diffeomorphisms of tori of dimension higher
than $1$. In the one-dimensional case, the celebrated
Denjoy theorem and examples establish a break in
dynamical behaviour at the regularity threshold $C^{1+
\mathrm{BV}}(\T^{1})$.\footnote{By $C^{1+\mathrm{BV}}
(\T^{1})$ we denote the space of circle diffeomorphisms
of the circle whose first derivative has bounded variation.}
A circle diffeomorphism with irrational
rotation number and regularity lower than
$C^{1+\mathrm{BV}}$ may have wandering intervals, while
a diffeomorphism of regularity $C^{1+\mathrm{BV}}$ cannot
(see, e.g. \cite{KatHass}).

In the one-dimensional case (see e.g. \cite{deMvanStr}
or \cite{Herm79}), a homeomorphism is assigned
a unique rotation number, and, as soon as it is irrational,
a continuous semi-conjugation to the rigid rotation can be
readily constructed. Denjoy's theorem is a rigidity theorem,
stating that if the homeomorphism is sufficiently
regular, the semi-conjugacy is in fact automatically
a continuous conjugacy. Arnold's theorem and subsequently
the Herman-Yoccoz theory (\cite{Herm79} and \cite{YocAst}) is a
further rigidity result in this setting, under additional
regularity and arithmetical assumptions.

It is not known whether Denjoy's theory admits a reasonable
generalization when the dimension of the torus is higher
than one, but it is certainly not directly generalizable (due to
the fact that the Denjoy-Koksma inequality fails, see
\cite{YocAst}).
A homeomorphism of a higher dimensional torus does not, in
general, have a unique rotation vector (see again \cite{Herm79}),
and even if this is the case, the minimality of the
corresponding translation does not imply the existence of
a continuous semi-conjugacy to it (e.g. \cite{Furst61}).

It is conjectured in \cite{McS93}, however, that
for diffeomorphisms $f$ of $\T ^{d}$, with $d\geq 2$, who
satisfy the \textit{additional} assumption that there
exists $\phi :\T^{d} \ra \T^{d}$, continuous and surjective
and such that
\begeq
\phi \circ f = R_{\a} \circ \phi
\endeq
with $R_{\a}$ minimal, a similar break should appear. That
is, it is possible to
construct such diffeomorphisms with wandering domains as
long as the diffeomorphism is in, say $C^{d+1-\e}$,
but not if it is more regular than, say,
$C^{d+1}$.

To our best knowledge, some examples of particular nature
and of regularity lower than the conjectured threshold have
been constructed by McSwiggen (see \cite{McS93} and \cite{McS95})
and by Sambarino and Passeggi (see \cite{PassegSamb13}).
McSwiggen's examples on $\T^{2}$ are based on
a Derived Anosov technique. A linear Anosov diffeomorphism of
$\T ^{3} $ is deformed in $\sm $ in order to turn
the saddle around the fixed point into a repeller. A diffeomorphism of
$\T ^{2}$ is then constructed
as the holonomy map along the unstable foliation that is
proved to survive the deformation. The radical loss of regularity
from $\sm$ to $C^{3-\e }$ comes from an inequality that has
to be satisfied by an algebraic function of
the eigenvalues of the original Anosov system and the need for
contraction in the functional space of some bundle sections,
and therefore in a quite
indirect manner. The diffeomorphism thus constructed is proved
to have wandering domains, while it is semi-conjugate to
the unstable holonomy map of the original Anosov
diffeomorphism by collapsing the repelling basin of
the origin to a point.
Moreover, even though it is not actually
mentioned in the paper, the rotation
vector is in fact Diophantine, since both the direction
and the modulus of the translation vector are given
by algebraic functions.

We think that, when compared to our result, this construction
represents an instance of the fact that in low regularity arithmetic
properties are irrelevant, while they become crucially so
above some finite (and hopefully universal and explicit) threshold.

\textbf{Acknowledgement}: This work was supported by the
ERC AdG grant no 339523 RGDD. The author would like to thank
Sebastian van Strien and Abed Bounemoura
for some useful discussions during the preparation of the
article.

\section{Notation}
\subsection{General notation}
By $F : \R ^{d}\ra \R ^{d}$ we denote a lift of $f$
(and in general the corresponding capital letter
will denote a lift whenever the small one denotes
a diffeomorphism of the torus). By  a tilde we denote the lift
of a point $x \in \T ^{d} \equiv \R ^{d} / \Z ^{d} $ to a representative in the
covering space $\tilde{x} \in \R ^{d}$.

A special case of diffeomorphisms of the torus is that of
translations. For $\a \in \T^{d}$, we define
\begeq
R_{\a} : x \mapsto x+\a \mod \Z ^{d}
\endeq

We will denote the space of $C^{s}$-smooth diffeomorphisms
that are isotopic to the identity
by $\mathrm{Diff}^{s}_{0}(\T ^{d})$, and the distance in
$\mathrm{Diff}^{s}_{0}$ between two diffeomorphisms
$f$ and $g$ by
\begeq
d_{s} (f,g) = \max _{0\leq \s \leq s} \| D^{\s} F - D^{\s}  G \|
 _{L^{\infty}}
\endeq
The space of $\sm$ diffeomorphisms will be furnished
with the corresponding topology.

If $\f : \T ^{d} \ra \R $, $\hat{\f}(k), k \in \Z ^{d}$ are
its Fourier coefficients, and
$N \in \N ^{*}$, we denote by
\begin{eqnarray*}
T_{N}\f (\. )&=& \sum _{|k| \leq N} \hat{\f}(k)e^{2i\pi \.} \\
\dot{T}_{N}\f (\. )&=& \sum _{0<|k| \leq N} \hat{\f}(k)e^{2i\pi \.} \\
R_{N}\f (\. )&=& \sum _{|k| > N} \hat{\f}(k)e^{2i\pi \.} 
\end{eqnarray*}
the inhomogeneous and homogeneous truncations and the
rest, respectively,
where $\Z ^{d}$ is equipped with the $\ell ^{1}$ norm.
The estimates
\begin{eqnarray*}
\| T_{N}\f (\. ) \| _{s} &\leq & C_{s} N^{s+d/2} \| \f \| _{0} \\
\| \dot{T}_{N}\f (\. ) \| _{s} &\leq & C_{s} N^{s+d/2} \| \f \| _{0} \\
\| R_{N}\f (\. ) \| _{s} &\leq & C_{s,s'} N^{-s' + s+d} \| \f \| _{s'},
\end{eqnarray*}
are well known, where $0\leq s \leq s'$.

If $f,g,u \in \mathrm{Diff}^{\infty}_{0}(\T ^{d})$, then,
see \cite{KrikAst},
\begin{eqnarray} \label{eq est comp}
\| g \circ f \| _{s} &\leq & C_{s} \| g \| _{s} (1+\| f \| _{s})
(1+ \| f \| _{0})^{s}
\\
\left\Vert g \circ (f+u)-\psi \circ f\right\Vert _{s} &\leq &
C_{s}\left\Vert gñ \right\Vert _{s+1}(1+\left\Vert f\right\Vert
_{0})^{s}(1+\left\Vert f\right\Vert _{s})\left\Vert u\right\Vert _{s}
\end{eqnarray}

Finally, the vector $ \a \in \T ^{d} $ is said to satisfy a Diophantine
condition of type $\gamma , \t $, if the following holds:
\begeq
\a \in DC(\gamma ,\t ) \Leftrightarrow
\forall k \in \Z ^{d} \setminus \{ 0\} ,
|k \. \a |_{\Z ^{d}} \geq \frac{\gamma ^{-1}}{|k|^{-1}}
\endeq

\subsection{Rotation vectors and sets} \label{subs rot vec}
For this paragraph, see \cite{MisZiem91} or \cite{Franks95}.
If $f \in \mathrm{Homeo}_{0}(\T ^{d}) \equiv \mathrm{Diff}_{0}^{0}(\T ^{d})$
and $x \in \T ^{d}$, we define $\r (x,f)$ as the following
limit, provided that it exists:
\begeq
\frac{F^{n}(\tilde{x})-\tilde{x}}{n}
\endeq
It is defined $\mod \Z ^{d}$, due to the arbitrary choice of a lift
for $f$.

We also define $\r (f)$, the rotation set of $f$, as
the accumulation points of
\begeq
\frac{F^{n_{i}}(\tilde{x_{i}})-\tilde{x_{i}}}{n_{i}}
\endeq
where $n_{i} \ra \infty$ and the $x_{i}  \in \T ^{d}$.
It can be shown that $\r (f)$ is the convex hull of
$\cup _{x \in \T ^{d}} \r (x,f)$.

If $\nu \in \mathcal{M}(f)$ (i.e. a probability measure on
$\T ^{d}$, invariant under $f$) then the quantity
\begeq
\int _{\T ^{d}} (F(\tilde{x}) - \tilde{x}) d\nu
\endeq
is well defined and denoted by $\r (\nu , f)$. The set
$\cup  \r (\nu ,f)$, where the union is over $\mathcal{M}
(f)$, is denoted by $\r _{meas}(f)$.

M. Herman defines the rotation set of a homeomorphism
precisely as $\r _{meas}(f)$ (his notation is different),
and his conditions on $f$ and the volume that it preserves
are needed in order to assure that
\begeq
\r _{meas}(f) \equiv \r (\nu , f)\equiv \r (\mu , f) \equiv \r (x,f)
\endeq
for every invariant measure $\nu$, and for every point $x\in \T ^{d}$,
as is automatically the case in the circle (see
\cite{Herm79}).

The condition imposed in thm. \ref{thm A} can be written
in the form\footnote{By $\mathrm{Conv}$ we denote the convex hull
of a set, i.e. the smallest closed convex set containing the
given one.}
\begeq
\a \in \mathrm{Conv}(\r (f)) = \r _{meas}(f)
\endeq
When $d =2$, it can be shown that $\r (f)$ is convex, see
\cite{MisZiem91}, so that $\r (f)  = \r _{meas}(f) $.

\section{Statement of the theorem}
We can now restate our main theorem in a more precise way.
\begin{theorem} \label{thm rig}
Let $d \in \N ^{*}$, $\gamma >0 $ and $\tau > d$. Then,
there exist $\e >0$ and $s_{0} >0$ such that if $\a \in \T ^{d}$
and $f \in \mathrm{Diff}^{\infty}_{0}(\T ^{d})$ satisfy
\begin{enumerate}
\item $\a \in DC(\gamma , \tau )$
\item \label{smallness cond item} $d_{0}(f(\. ) , R_{\a} ) < \e$ and
$d_{s_{0} } (f(\. ) , R_{\a} ) < 1$
\item  \label{rot set item}$\a \in \r _{meas}(f) $
\end{enumerate}
then $f $ is $\sm$ conjugate to $R_{\a}$. Moreover, the
conjugation can be chosen close to the $Id$.
\end{theorem}
Since such results tend to generalize to finite differentiability,
we expect the following conjectural theorem to be true.
\begin{theorem}[Conjectural] \label{thm rig finite dif}
Let $d \in \N ^{*}$, $\gamma >0 $ and $\tau > d$. Then,
there exist $\e >0$ and $\kappa , s_{0} >0$ with
$0<\kappa  <s_{0}$ such that, if $\a $
and $f \in \mathrm{Diff}^{s}(\T ^{d})$ with $s>s_{0}$ satisfy the conditions of theorem \ref{thm rig} in items
$1-3$,
then $f$ is $C^{s-\kappa}$ conjugate to $R_{\a}$.
The conjugation can be chosen close to the $Id$.
\end{theorem}
This conjectural theorem is implied, for instance, by  the proof in \cite{Herm79}, which
is carried out by approximation of finitely differentiable mappings by analytic ones.
The proof is valid as we point out in the answer to the following question.

\begin{question}
Does thm \ref{thm rig} hold true in the
real analytic category?
\end{question}

The is of course yes and an easy but non-optimal argument is as follows. If the
perturbation is small enough in some analytic norm, then it
is small enough in $\sm$. Therefore, thm \ref{thm rig}
applies and the diffeomorphism is $\sm$ conjugate to the
rigid rotation $R_{\a}$. As a consequence, its rotation set
is reduced to $\{ \a \} $,
and M. Herman's proof can be applied by just dropping the
volume preservation assumption.

\section{Proof of theorem \ref{thm rig}}

The proof relies on two lemmas. The first one is a K.A.M.
lemma of very classical flavour and estimates, and represents
one step of the K.A.M. scheme that constructs successive
conjugations reducing the diffeomorphism $f$ to the rigid
rotation $R_{\a}$. The second lemma is a geometric one and
relates the size of the perturbation with $\r _{meas}$.

For the scheme to
produce a converging product of conjugations, two conditions
are needed. The first, more standard one, is a closeness to a
rotation condition in an appropriate topology, and its general
form is like the one of item \ref{smallness cond item} of the statement of thm \ref{thm rig}.
The second one is used in making sure that the perturbed
diffemorphism $f$ does not drift away from $R_{\a}$: clearly,
if $\beta$ is a vector with rational coordinates, very close
to $\a$, the two corresponding rotations are not conjugate.\footnote{
We remark that the rotation set of a diffeomorphism of $\T ^{d}$
is only preserved by conjugations that are isotopic to the $Id$.}
Such a condition can be imposed on the rotation set of $f$,
the perturbed diffeomorphism, and a possible condition
would then be
\begeq
\r (f) \equiv \{ \a \}
\endeq
where $\r (f) $ is defined in paragraph \ref{subs rot vec}.
In fact, the following weaker condition would be sufficient:
\begeq
\exists x \in \T ^{d}, \r (x,f) = \{ \a \}
\endeq
However, it is enough that $\a \in \mathrm{Conv} \r (f) = \r _{meas}(f)$, which is
exactly condition in item \ref{rot set item} of the theorem.

All three conditions are weaker than
\begeq
\exists \Phi :\R ^{d} \ra \R ^{d}, \| \Phi - Id \|
_{L^{\infty}} < \infty
\endeq
a measurable mapping, such that
\begeq
\Phi \circ F = R _{\a} \circ \Phi
\endeq
This last condition implies that for $a.e.$ $x \in \T ^{d}$
\begeq
\r (x,f) = \{ \a \}
\endeq
If $\Phi$ is assumed to be continuous, this holds for every
$x \in \T ^{d}$ and in fact $\r (f) = \{\a \}$.

Finally, this last condition is weaker than the existence of
$ \phi :\T ^{d} \ra \T ^{d}$, surjective and respectively
measurable or continuous, such that
\begeq
\phi \circ f = \R _{\a} \circ \phi
\endeq
In the one-dimensional case, the existence of such a
semi-conjugation is automatic as soon as the rotation number
of the homeomorphism is irrational. In the higher-dimesnional
case, however, the existence of a semi-conjugation to a
minimal rotation is an additional and restrictive hypothesis, and
brings us back to the context of the conjecture mentioned in the
introduction.

\subsection{Inductive lemma}
We now state and recall the proof of the inductive conjugation lemma.
\begin{lemma} \label{ind KAM lem}
Let $\a \in DC(\gamma , \tau ) \subset \T ^{d}$ and
$f \in \mathrm{Diff}^{\infty }(\T ^{d})$, and call
$\| f - R_{\a} \|_{C^{s}} = \e _{s}$.
Then, for some absolute constant $C>0$ and for every $N \in \N ^{*}$ such that
\begeq
C \gamma N^{2\t +d+2 }\e _{0}<1
\endeq
there exists $\phi \in \mathrm{Diff}^{\infty}(\T ^{d})$ such that 
\begeq
\phi \circ f \circ \phi ^{-1} = f'
\endeq
and the following hold true for the diffeomorphism $f' \in
\mathrm{Diff}^{\infty} $ thus defined.

There exists $\beta \in \T ^{d}$ such that
$\| f' - R_{\beta} \|_{C^{s}} = \e ' _{s}$ satisfies
\begeq
\e ' _{s} \leq C_{s,s'}\left(
N^{s+2\t +d+2}\e _{0}^{2}+
N^{\t +d/2}\e _{0}\e_{s}+
N^{s-s'+d} (1+ N^{s+ \t +d/2}\e_{0} )\e _{s'}
\right)
\endeq
for every $0\leq s \leq s' < \infty$.

Moreover, the conjugation $\phi$ satisfies
\begeq
\| \phi \|_{s} \leq C_{s}\gamma N^{s+\t + d/2}\e _{0}
\endeq
\end{lemma}

Naturally, $\beta \simeq \a + \int (f - R_{\a } )d\mu$ and it represents a drift
of the perturbed diffeomorphism with respect to  $R_{\a}$.
The constants appearing in the statement depend on
$ \t $ and $ d$, but not on $N$.
The proof is classical, but we sketch it for the sake of completeness.
\begin{proof}
Let $d=2$ in order to simplify notation, without any loss of generality. Then,
\begeq
f (\. ) = R_{\a} +
\begin{pmatrix}
f _{1} (\. ) \\
f _{2} (\. )
\end{pmatrix}
\endeq
where $f _{i} (\. ) : \T ^{2} \ra \R $ for $i = 1,2$,
are small in the $\sm$ topology.

If we call
\begeq
\phi (\. ) = Id +
\begin{pmatrix}
\phi _{1} (\. ) \\
\phi _{2} (\. )
\end{pmatrix}
\endeq
where $\phi _{i} (\. ) : \T ^{2} \ra \R $ for $i = 1,2$,
then, for the conclusion of the lemma to be true,
they need only satisfy the equation
\begeq
\phi _{i} (\. ) \circ R_{\a} - \phi _{i} (\. )
+ \dot{T}_{N} f_{i} (\. ) = 0
\endeq
Such functions $\phi _{i}$ exist and are uniquely defined
in $\sm _{0} (\T ^{2} )$. They satisfy the estimate
\begeq
\| \phi _{i} (\. ) \| _{s} \leq C_{s} \gamma
N^{s+\tau + d/2} \| f_{i}\| _{0}
\endeq
Then, we can calculate
\begin{eqnarray*}
\phi \circ f \circ \phi ^{-1} &=&
\left(
Id +
\begin{pmatrix}
\phi _{1} (\. ) \\
\phi _{2} (\. )
\end{pmatrix}
\right) 
\circ
\left(
R_{\a} +
\begin{pmatrix}
f _{1} (\. ) \\
f _{2} (\. )
\end{pmatrix}
\right) \circ
\left(
Id -
\begin{pmatrix}
\phi _{1} (\. ) \\
\phi _{2} (\. )
\end{pmatrix}
+ O(\phi _{i}^{2})
\right) \\
&=& \left(
Id +
\begin{pmatrix}
\phi _{1} (\. ) \\
\phi _{2} (\. )
\end{pmatrix}
\right) 
\circ
\left(
R_{\a}-
\begin{pmatrix}
\phi _{1} (\. ) \\
\phi _{2} (\. )
\end{pmatrix}
+ T_{N}
\begin{pmatrix}
f _{1} (\. ) \\
f _{2} (\. )
\end{pmatrix} 
+ O(\. )
\right) \\
&=&
R_{\a} +
\begin{pmatrix}
\hat{f} _{1} (0) \\
\hat{f} _{2} (0)
\end{pmatrix} +
\begin{pmatrix}
\phi _{1} (\. ) \\
\phi _{2} (\. )
\end{pmatrix}
\circ R_{\a}
-
\begin{pmatrix}
\phi _{1} (\. ) \\
\phi _{2} (\. )
\end{pmatrix}
+ \dot{T}_{N}
\begin{pmatrix}
f _{1} (\. ) \\
f _{2} (\. )
\end{pmatrix} 
+ O(\. )  \\
&=& 
R_{\a} +
\begin{pmatrix}
\hat{f} _{1} (0) \\
\hat{f} _{2} (0)
\end{pmatrix} + O(\phi _{i}^{2}, \partial (\dot{T}_{N} f _{i})
 . \phi _{i}, \partial (\phi _{i}) . \phi _{i},
R_{N} f _{i} \circ (Id - \phi _{i}))
\end{eqnarray*}
The $O(\. )$ term in the last line, which, anticipating the
next section we call
\begeq
\begin{pmatrix}
f ^{\prime} _{1} (\. ) \\
f ^{\prime} _{2} (\. )
\end{pmatrix} 
\endeq
so that
\begeq
\phi \circ f \circ \phi ^{-1} =
R_{\a} +
\begin{pmatrix}
\hat{f} _{1} (0) \\
\hat{f} _{2} (0)
\end{pmatrix}
+
\begin{pmatrix}
f ^{\prime} _{1} (\. ) \\
f ^{\prime} _{2} (\. )
\end{pmatrix}
\endeq
can be estimated in the $C^{s}$-norm by
\begeq
C_{s,s'}\left(
N^{s+2\t +d+2}\e _{0}^{2}+
N^{\t +d/2}\e _{0}\e_{s}+
N^{s-s'+d} (1+ N^{s+ \t +d/2}\e_{0} )\e _{s'}
\right)
\endeq
This concludes the proof of the lemma.
\end{proof}

\subsection{A posteriori estimate on the obstruction}

The following elementary and well known observation
establishes a relation between the displacement of points in
the torus $\T ^{d}$ under a diffeomorphism $g$ with its
rotation set $\r (g)$.

\begin{lemma} \label{conv lem}
Let $g \in \mathrm{Diff}^{\infty}_{0} (\T ^{d})$ and $\beta \in \T ^{d} $.
If there exists $x \in \T ^{d}$ such that $ \r (x,g) = \{ \beta \}$,
then
$\beta \in \mathrm{Conv} \{ G (\tilde{x}) -\tilde{x} ,x \in \T^{d} \}$.
\end{lemma}
Inspection of the proof shows that the condition on the existence
of an orbit rotating at speed $\beta $ can be relaxed to $\beta \in \r (g)$.

\begin{proof}
Let $x \in \T ^{d}$ be such that $ \r (x,g) = \{ \beta \}$. Then,
\begeq
\frac{G^{n} (\tilde{x}) -\tilde{x}}{n} =
\frac{G (G^{n-1} (\tilde{x})) - G^{n-1} (\tilde{x})+ G (G^{n-2} (\tilde{x})) -
G^{n-2} (\tilde{x}) + \cdots + G( \tilde{x}) - \tilde{x}}{n}
\endeq
converges to $\beta$.
Since the right-hand side is an element of
$\mathrm{Conv} \{ G (\tilde{x}) -\tilde{x} ,x \in \T^{d} \}$ and the latter
set is closed, the lemma is proved.
\end{proof}

In the context of lemma \ref{ind KAM lem}, we obtain the
following corollary.
\begin{corollary} \label{a post lem}
There exists an absolute constant $C>0$ depending only on $d$ such that,
under the hypotheses of lem. \ref{ind KAM lem}, and assuming
additionally that $\a \in \mathrm{Conv} ( \r (f) )$, 
\begeq
\left\lVert
\begin{pmatrix}
\hat{f} _{1} (0) \\
\hat{f} _{2} (0)
\end{pmatrix}
\mod \Z ^{2}
\right\rVert
\leq C \e '_{0}
\endeq
\end{corollary}

In the proof we assume for simplicity that $ \r (x,f) = \{ \a \}$ for
some $x \in \T ^{2}$.
The proof of the corollary as it is stated follows easily.
\begin{proof}
If $ \r (x,f) = \{ \a \}$, then $\r (\f (x),f') = \{ \a \}$. Then, by
lemma \ref{conv lem},
\begeq
\a \in \mathrm{Conv}(F'(\. ) - Id)
\endeq
Consequently,
\begeq
0 \in \mathrm{Conv}(F'(\. ) - \a ) =
\mathrm{Conv}
\left(
\begin{pmatrix}
\hat{f} _{1} (0) \\
\hat{f} _{2} (0)
\end{pmatrix}
+
\begin{pmatrix}
f ^{\prime} _{1} (\. ) \\
f ^{\prime} _{2} (\. )
\end{pmatrix}
\right)
\endeq
The corollary follows directly.
\end{proof}

Since invariant measures are accumulation points of Dirac
measures uniformly distributed on finite segments of orbits,
we can relax the condition $\a \in \r (f) $ to
$\a \in \r _{meas}(f) $: We must have
\begeq
0 = \int _{\T ^{d}}
\left(
\begin{pmatrix}
\hat{f} _{1} (0) \\
\hat{f} _{2} (0)
\end{pmatrix}
+
\begin{pmatrix}
f ^{\prime} _{1} (\. ) \\
f ^{\prime} _{2} (\. )
\end{pmatrix}
\right)
d(\phi_{*}\nu )
\endeq
for every (fixed) $\nu \in \mathcal{M}(f)$ such that
$\r (\nu ,f)= \a $. We immediately get the same estimate as
in cor. \ref{a post lem}.

\subsection{KAM scheme and convergence}

The estimates provided by lemma \ref{ind KAM lem} are
sufficient for the convergence of the corresponding scheme,
provided that some smallness conditions are satisfied,
and that we can linearize around the same
rotation $\a$ throughout the scheme, so that no
"counter-term" is needed (as in the normal form version of the
theorem in \cite{Herm79}), and the Diophantine condition can
be kept constant throughout the scheme. This second condition
is assured by the existence of an orbit rotating like
$R_{\a}$, and by corollary \ref{a post lem}.

Let us state this formally in the following proposition.
\begin{proposition} \label{prop conv KAM}
Let $\a \in DC(\gamma , \tau ) \subset \T ^{d}$ and
$f = f_{1} \in \mathrm{Diff}^{\infty }(\T ^{d})$, and call
$\| f_{1} - R_{\a} \|_{C^{s}} = \e _{s,1}$.
Then, there exist $\epsilon >0 $ and $s_{0} \in \N^{*} $
such that if
\begeq
\e _{0,1}<\epsilon \text{ and } \e _{s_{0},1} <1
\endeq
and if
\begeq
\a \in \r (f_{1})
\endeq
then there exist inductively defined sequences $\phi _{n} \in
\mathrm{Diff}^{\infty }_{0}(\T ^{d})$
and $f _{n} \in \mathrm{Diff}^{\infty }_{0}(\T ^{d})$ such that 
\begeq
\phi _{n} \circ f _{n} \circ \phi ^{-1}_{n} = f_{n+1}
\endeq
with
\begin{equation} \label{eq cvgence}
\| f_{n} - R_{\a} \|_{C^{s}} = \e _{s,n}  \xrightarrow[n\ra \infty]{} 0 , \forall s \in \N
\end{equation}
Moreover,
\begeq
\prod _{k=1}^{n} \phi _{k}  \xrightarrow[n\ra \infty]{} \phi
\in \mathrm{Diff}^{\infty }(\T ^{d})
\endeq
\end{proposition}
Clearly, this proposition implies thm. \ref{thm rig}. The
proposition is proved by iteratively applying lemma
\ref{ind KAM lem} and then corollary \ref{a post lem}
in the following, now classical, way.

\begin{proof}
Let $N = N_{1} \in \N ^{*}$ be large enough, chose
$\s >0$, and define inductively
\begeq
N_{n} = N_{n-1}^{1+\s} = N^{(1+\s)^{n-1}}
\endeq
to be the order of truncation at the $n$-th step, as in the proof of
lem. \ref{ind KAM lem}.

Assume, now, that $\phi _{n-1}$ has already been
constructed, so that $f_{n}$ is well defined. Suppose,
additionally, that $f_{n}$ satisfies the hypotheses
of lem. \ref{ind KAM lem} for $\e _{s} = \e _{s,n}$
and $N=N_{n}$. Then, application of lem. \ref{ind KAM lem}
and then of cor. \ref{a post lem} grants the existence of
$\phi _{n}$ such that
\begeq
f_{n+1} = \phi _{n} \circ f _{n} \circ \phi _{n}^{-1}
\endeq
satisfies, for all $0 \leq s \leq s' < \infty $, the inequality
\begin{dmath} \label{ind ineq}
\e _{s,n+1} \leq C_{s,s'}\left(
N_{n}^{s+2\t +d+2}\e _{0,n}^{2}+
N_{n}^{\t +d/2}\e _{0,n}\e_{s,n}+
N_{n}^{s-s'+d} (1+ N_{n}^{s+\t +d/2}\e_{0,n} )\e _{s',n}
\right)
\end{dmath}
We remind that $\e _{n,s} $, defined after eq. \ref{eq cvgence},
represents the distance of the diffeomorphism $f_{n}$ from
the \textit{fixed} rotation $R_{\a}$, thanks to the a posteriori
estimate of cor. \ref{a post lem}.

Since the term in $\e_{0,n} \e _{s',n}$ is not present in thm. $10$
of \cite{FK2009}, we partially reprove the convergence of the
scheme. In fact, only the two main steps, lemmata $11$ and
$14$ of the reference, have to
be proved for the kind of estimates that we have herein\footnote{
The discrepancy is due to the fact that in the reference
the dynamical system considered is a cocycle, and the last
term herein comes from composition of mappings, whereas in
the context of cocycles mappings are composed uniquely with
the $\exp$ and only products need to be considered.}.
We therefore need the following lemma, which is proved in the
appendix, section \ref{appendix}.
\begin{lemma} \label{lem KAM 1}
Let $\e _{s,n}$ satisfy the inductive estimates of eq. \ref{ind ineq}. If,
moreover
\begin{eqnarray*}
\e _{0,1} &<& N_{1}^{-\g _{0} } \\
\e _{s_{0},1} &<& N_{1}^{b}
\end{eqnarray*}
for some appropriately chosen $\g _{0},b>0$, then the double
sequence $\e _{s,n}$ is well defined and for all $n$
\begin{eqnarray*}
\e _{0,n} &<& N_{n}^{-\g _{0}} \\
\e _{s_{0},n} &<& N_{n}^{b}
\end{eqnarray*}
\end{lemma}
We note, en passant, that this lemma implies that, under the relevant
smallness conditions on $\e _{s,0}$, the smallness conditions of
lemma \ref{ind KAM lem} are satisfied for all $n$. Therefore, the
double sequence $\e _{s,n}$ is well defined and we only need to
establish its convergence.

We then show that, given the decay and growth
rates granted by the previous lemma, we can actually do slightly better.
The following lemma allows us to bootstrap exactly like in \cite{FK2009}
and conclude the convergence.
\begin{lemma} \label{lem KAM 2}
Let $\e _{s,n}$ satisfy the inductive estimates of eq. \ref{ind ineq}. If,
moreover
\begin{eqnarray*}
\e _{0,n} &=& O(N_{n}^{-\g _{0}}) \\
\e _{s_{0},n} &=& O(N_{n}^{b})
\end{eqnarray*}
with $\g _{0}$, $b$ and $s_{0} $ as in the previous lemma
then, there exist $\w _{0} , \w >0$ such that
\begin{eqnarray*}
\e _{0,n} &=& O(N_{n}^{-(1+\w _{0} )\g _{0}}) \\
\e _{(1+\w) s_{0},n} &=& O(N_{n}^{b})
\end{eqnarray*}
\end{lemma}

Thus, we have successively proved the following assertions:
\begin{itemize}
\item at each step, the inductive smallness hypothesis
for lem. \ref{ind KAM lem} to be applicable is satisfied by
$f_{n}$ and $N_{n}$. Therefore, the double sequence
$\e _{s,n}$ is well defined for all $n\in \N^{*}$ and
$s\in \N$.
\item for every fixed $s\in \N$ and $\l \geq 0$
\begeq
N^{\l}_{n}\e _{s,n}   \xrightarrow[n\ra \infty]{} 0
\endeq
The shorthand $f _{n} - R_{\a} = O _{C^{\infty}} (N_{n}^{-\infty})$
is common. This is obtained by the fast convergence of $\e _{0,n}$ to $0$,
faster than any power of $N_{n}$, and the slow growth of $\e _{s,n}$
(as a fixed power of $N_{n}$) for every $s$ fixed. Convexity estimates
allow us to conclude that for every $0<s'<s$, $\e _{s',n}\ra 0$.
\item the fast convergence of $\e _{s,n}$ to $0$ and the fact
that
\begeq
\| \Phi _{n} - Id \|_{s} \leq C_{s} \gamma N_{n}^{s+\t + d/2}\e _{0,n}
\endeq
imply, together with eq. \ref{eq est comp} that the product of
successive conjugations
\begeq
\prod _{k=1}^{n} \phi _{k} \in  \mathrm{Diff}^{\infty} (\T^{d})
\endeq
converges in the $\sm $ topology to a well defined diffeomorphism
$\phi $.
\end{itemize}
This concludes the proof of the proposition, and thus of thm \ref{thm rig}.
\end{proof}

\section{A remark on the proof}
In the one-dimensional case, the theory of the rotation
number for an orientation preserving homeomorphism of
$\T ^{1}$ is considerably stronger, thanks to the existence
of a natural cyclic order (or of a total order in the
covering space $\R ^{1}$). The analogue of our argument in
the one-dimensional case would be the following. Consider
a Diophantine rotation $\a$ and perturb it. Suppose that one
orbit with rotation number $\a$ survives under perturbation.
Solve the linear cohomological equation and observe that if
the obstruction (the mean value of the perturbation) is not
of "second order", a contradiction would be established,
e.g. by fitting a rational number between
\begeq
\mathrm{Conv} \{ F (\tilde{x}) -\tilde{x} ,x \in \T^{1} \}=
 \{ F (\tilde{x}) -\tilde{x} ,x \in \T^{1} \}
\endeq
and $\a$, see \cite{KatHass}.

Of course, and this is a particularity of the one-dimensional
theory, if one orbit has rotation number $\a$, then all
orbits do. This, however, is not an essential part of the
proof of the existence of a smooth conjugacy to the rigid
rotation, since the proof of the uniqueness of the rotation
number is formally independent of
the construction of the K.A.M. scheme and of the proof
of its convergence.
M. Herman in his thesis, \cite{Herm79}, defines the rotation
number for circle diffeomorphisms using an invariant measure
instead of a combinatorial definition as in \cite{KatHass} or \cite{deMvanStr}.
Accordingly to the one-dimensional case, his definition of the rotation set
for diffeomorphisms preserving a volume form assures that the rotation set is
thus reduced to a point, and this hypothesis is, in fact,
needlessly strong.

\section{Appendix} \label{appendix}

We now provide the missing proofs of lemmata \ref{lem KAM 1} and
\ref{lem KAM 2}. We note that, when eq. \ref{ind ineq} is
compared with the corresponding eq. $7.2$ of
\cite{FK2009}, which in our notation reads
\begin{dmath}
\e _{s,n+1} \leq C_{s,s'}\left(
N_{n}^{a+Ms}\e _{0,n}^{1+\s _{0}}+
N_{n}^{a+ms}\e _{0,n}\e_{s,n}+
N_{n}^{a - (s'-s)\mu } (\e _{s',n}+ N_{n}^{\bar{\mu} s'}\e_{0,n} )
\right)
\end{dmath}
the agreeing terms correspond to the
admissible choice of parameters (in the notation of
the reference)
\begin{eqnarray*}
\s _{0} =  1 &,& a = 2\t +d+2 \\
M = 1 &,& m = 0 \\
\mu = 2 &,& \bar{\mu } = 0
\end{eqnarray*}
However, there do not seem to exist admissible values of the parameters
$g $, $\s $ and $\kappa$ of the reference for which either our estimates can
be brought to the form of those of the reference, or for which the proof found
therein produces convergence of the scheme for our type of estimates.
Consequently, we take up the
proof of convergence for our type of estimates and we remark that
the broader values of parameters (or the broader scope of types of estimates)
are obtained thanks to the fact that we consider $a$ not as an "affine" parameter,
but as a "homogeneous" one, when compared to $\g _{0}$, $b$ and $s_{0}$.
%

The additional term appearing in eq. \ref{ind ineq}, namely
$N_{n}^{2s-s'+a}\e_{0,n} \e _{s',n}$,
is due to the composition of the conjugation with
the perturbation (and is not present in the context of
cocycles). By defining $\g _{0}$, $b$ and $s_{0}$ as multiples of $a$ we manage
to absorb the additional growth by $N_{n}^{s}$ into $\e _{0,n}$, which
is shown to decay fast enough.

Let us now proceed to the actual proofs of the lemmata.

\begin{proof}[Proof of lem. \ref{lem KAM 1}]
In view of the estimate of eq. \ref{ind ineq} by setting $s=0, s' =s_{0}$
and $s=s'=s_{0}$, we only need to show that if $\e_{0,n}<N_{n}^{-\g _{0}}$ and
$\e_{s_{0},n}<N_{n}^{b}$, then
\begin{dmath*}%
\e _{0,n+1} \leq C_{0,s_{0}}\left(
N_{n}^{a}\e _{0,n}^{2}+
N_{n}^{a/2}\e _{0,n}^{2}+
N_{n}^{-s_{0}+a/2} (1+ N_{n}^{a/2}\e_{0,n} )\e _{s_{0},n}
\right)
\end{dmath*}
and
\begin{dmath*}%
\e _{s_{0},n+1} \leq C_{s_{0},s_{0}}\left(
N_{n}^{s_{0}+a}\e _{0,n}^{2}+
N_{n}^{a/2}\e _{0,n}\e_{s_{0},n}+
N_{n}^{a/2} (1+ N_{n}^{s_{0}+a/2}\e_{0,n} )\e _{s_{0},n}
\right)
\end{dmath*}
satisfy $\e_{0,n+1}<N_{n+1}^{-\g _{0}}$ and $\e_{s_{0},n+1}<N_{n+1}^{b}$.

The inequalities that need to be verified in the limit of $N_{1}$ large enough read
\begin{eqnarray*}
a- 2\g _{0} < -(1+\s ) \g _{0} &,& a/2 - 2\g _{0} < -(1+\s ) \g _{0} \\
-s_{0}+a/2 + b < -(1+\s ) \g _{0} &,& -s_{0}+a - \g _{0} + b < -(1+\s ) \g _{0}  \\
s_{0}+a - 2\g _{0} < (1+\s )b &,& a/2 - \g _{0} + b < (1+\s )b \\
a/2 +b < (1+\s )b &,& s_{0} +a - \g _{0} +b < (1+\s )b
\end{eqnarray*}%
Substitution by $\g _{0} = \l a$, $s_{0} = \mu a$ and $b = \nu a$ where
$\l , \mu , \nu$ are positive reals, gives%
\begin{subequations}
\begin{eqnarray}
(1-\s ) \l &>&1 \label{ineq KAM Ia} \\
\mu - (1+\s ) \l - \nu &>& 1/2  \label{ineq KAM Ib} \\
\mu - \nu -\s \l  &>&1 \label{ineq KAM Ic}\\
2\l + (1+\s )\nu -\mu  &>&1 \label{ineq KAM Id} \\
\s \nu + \l  &>&1 \label{ineq KAM Ie} \\
\s \nu  &>& 1/2 \label{ineq KAM If} \\
\mu - \s \nu - \l  &>& 1 \label{ineq KAM Ig}
\end{eqnarray}
\end{subequations}
Inequality \ref{ineq KAM Ia} implies that $\l >1$, so that ineq. \ref{ineq KAM Ie}
is redundant. For ineq. \ref{ineq KAM Ib}, \ref{ineq KAM Ic}, \ref{ineq KAM Id} and \ref{ineq KAM Ig}
to be compatible, we need
\begin{eqnarray*}
(1+\s ) \l + \nu + 1/2 &<& 2\l + (1+\s )\nu -1 \\
1+ \s \l + \nu &<& 2\l + (1+\s )\nu -1  \\
1+\s \nu +\l &<& 2\l + (1+\s )\nu -1
\end{eqnarray*}
or equivalently
\begin{eqnarray*}
3/2 &<& (1-\s )\l +\s \nu \\
2 &<& (2-\s )\l + \s \nu   \\
2 &<& \l + \nu 
\end{eqnarray*}
The first two ineq. are implied by \ref{ineq KAM Ia}, \ref{ineq KAM Ie} and \ref{ineq KAM If}. We thus
need to impose
\begin{subequations}
\begin{eqnarray}
1 &>& \s \label{choice param a} \\
\l + \nu &>&2 \label{choice param b} \\
(1-\s ) \l &>&1\label{choice param c}\\
\s \nu  &>& 1/2 \label{choice param d}
\end{eqnarray}
\end{subequations}
and choose $\mu$ such that
\begeq
\max \{ (1+\s ) \l + \nu + 1/2 ,1+ \s \l + \nu , 1+\s \nu +\l \} < \mu <  2\l + (1+\s )\nu -1
\endeq
The conditions are equivalent to
\begin{eqnarray*}
\frac{\l -1}{\l} > \s &>& \frac{1}{2\nu}  \\
\l &>&\frac{2\nu}{2\nu -1} \\
\nu &>&1/2
\end{eqnarray*}
\end{proof}

At this point, we encourage the reader to check lemmata $12$ and $13$ in \cite{FK2009},
as they prepare the following proof of the lemma corresponding to lem. $14$ of the
reference.

\begin{proof}[Proof of lem. \ref{lem KAM 2}]
Let $\e _{0,n} < \bar{C} N_{n}^{-\g _{0}}$ and
$\e _{s_{0},n} < \bar{C} N_{n}^{b}$. We first prove that there exists $\w _{0}$ such that
\begeq
\e _{0,n} = O(N_{n}^{-(1+\w _{0})\g_{0}})
\endeq
We calculate directly
\begin{eqnarray*}
\e _{0,n+1} &\leq & \bar{C} C_{0,s_{0}}\left(
N_{n}^{a - \g _{0}}\e _{0,n}+
N_{n}^{-s_{0}+a/2 + b}  \right) \\
&\leq &  \bar{C} C_{0,s_{0}} ^{2} \left(
N_{n}^{a - \g _{0}}\left(
N_{n-1}^{a - \g _{0}}\e _{0,n-1}+
N_{n-1}^{-s_{0}+a/2 + b}  \right)+
N_{n}^{-s_{0}+a/2 + b}  \right) \\
&\leq & ( \bar{C} C_{0,s_{0}} )^{2} \left(
N_{n}^{a - \g _{0}}\left(
N_{n-1}^{a - 2\g _{0}}+
N_{n-1}^{-s_{0}+a/2 + b}  \right)+
N_{n}^{-s_{0}+a/2 + b}  \right) \\
&\leq & ( \bar{C} C_{0,s_{0}} )^{2} \left(
N_{n}^{\frac{2+\s}{1+\s } a - \frac{3+\s}{1+\s }\g _{0} } +
N_{n}^{\frac{3+2\s}{2+2\s }a - \g _{0} - \frac{s_{0} - b}{1+\s}} +
N_{n}^{-s_{0}+a/2 + b}  \right)
\end{eqnarray*}
Again in the limit $n \ra \infty$ we only need to verify that
\begin{eqnarray*}
\frac{2+\s}{1+\s } a - \frac{3+\s}{1+\s }\g _{0}  &<& -(1+\s )(1+\w _{0} )\g _{0}  \\
\frac{3+2\s}{2+2\s }a - \g _{0} - \frac{s_{0} - b}{1+\s}  &<& -(1+\s )(1+\w _{0} )\g _{0}  \\
-s_{0}+a/2 + b  &<& -(1+\s )(1+\w _{0} )\g _{0} 
\end{eqnarray*}

The first inequality holds true for
\begeq
0<\w _{0} < \frac{2\l - \s - 2 - \l \s (1+\s )}{\l (1+\s )^{2}}
\endeq

The second holds true as long as
\begeq
\mu > \frac{3}{2} + \s + \nu + (1+ \s )((1+\s )(1+\w _{0} )-1)\l
\endeq
This inequality is verified for $\w _{0} = 0$, by
\ref{choice param a} and \ref{choice param c} and by 
the choice of $\mu$, and therefore also verified for $\w _{0}$
small enough.

The third inequality is equivalent to
\begeq
\mu > \frac{1}{2} + \nu + (1+\s )(1+\w _{0})\l
\endeq
which is verified provided that $\w _{0}$ is small enough, by the choice of $\mu $.

The second assertion of the lemma follows directly from gain in the speed of convergence
for $\e _{0,n}$ (i.e. from the fact that we can replace $\g _{0}$ by $(1 + \w  _{0}) \g  _{0}$
for some fixed $\w  _{0} >0 )$), and the fact that in the inequalities $s_{0}$
\ref{ineq KAM Ib}, \ref{ineq KAM Ic},\ref{ineq KAM Id} and \ref{ineq KAM Ig} $\mu$ and
$\l$ appear with opposite signs, so that increase in $\mu $ can be compensated
by the increase in $\l \ra (1+\w _{0})\l $ without the inequality being violated for the same choice of $\nu$.
The eventual multiplicative constants can be absorbed in the exponents in the limit of
$n $ large enough.
\end{proof}
This last lemma can be used in order to boostrap and obtain prop. \ref{prop conv KAM}.
%
%


\bibliography{aomsample}
\bibliographystyle{aomalpha}

\end{document}